\documentclass[10pt]{article}
\font\smallit=cmti10

\usepackage{amssymb,latexsym,amsmath,epsfig,amsthm, geometry} 
\makeatletter
\numberwithin{equation}{section}
\renewcommand\section{\@startsection {section}{1}{\z@}
{-30pt \@plus -1ex \@minus -.2ex}
{2.3ex \@plus.2ex}
{\normalfont\normalsize\bfseries\boldmath}}

\renewcommand\subsection{\@startsection{subsection}{2}{\z@}
{-3.25ex\@plus -1ex \@minus -.2ex}
{1.5ex \@plus .2ex}
{\normalfont\normalsize\bfseries\boldmath}}

\renewcommand{\@seccntformat}[1]{\csname the#1\endcsname. }

\makeatother

\newtheorem{theorem}{Theorem}

\theoremstyle{definition}

\newtheorem{remark}{Remark}
\newcommand{\Mod}[1]{\ (\mathrm{mod}\ #1)}
\begin{document}

\begin{center}
\uppercase{\bf An Elementary Approach For Sums of Consecutive cubes being prime power squares}
\vskip 20pt
Atilla Akkuş
\\
{\smallit Ankara, Türkiye}\\
{\tt atilla.akkus@ug.bilkent.edu.tr}
\end{center}
\vskip 20pt
\vskip 30pt
\centerline{\bf Abstract}
\noindent
This paper proposes an elementary solution to a special case of finding all perfect squares that can be written as sum of consecutive integer cubes. It is shown that there are no non-trivial solutions if the perfect square is a prime power, i.e., it is divisible by two different primes if a non-trivial one exists. Solution mostly depends on $v_{p}(x)$ and general forms of Pythagorean triples.
\pagestyle{myheadings}
\thispagestyle{empty}
\baselineskip=12.875pt
\vskip 30pt
\section{Introduction}
A fundamental theorem in number theory states the equality between the square of $n$th triangular number ($\Delta_{n}$) and the sum of the first $n$ positive integer cubes. The general problem emerges when the theorem is extended to determine all cases when a perfect square equals to the sum of consecutive positive integer cubes of any interval. In other words, to find all positive\footnote{Notice that assuming $x,k,y$ positive will not cause loss of generality.} integer triples $(x,k,y)$ such that
\begin{equation}
    x^3 + (x+1)^3 + ... + (x+k-1)^3 = y^2. \label{maineq}
\end{equation}
\\
A remarkable special case ($k$ = 3) was solved by Uchiyama in 1977 \cite{Uchiyama}. Later, in 1985, Cassels proposed another solution for the same problem, depending on the elliptic curves \cite{Cassels}. The first holistic approach appears a decade after in Stroeker's paper \cite{Stroeker}, in which all solutions to Eq. (\ref{maineq}) under the condition $2\le k \le 49$ or $k = 98$ are given with the aid of computers. Complete solutions of the equation became achievable to find with Pletser's recursive method \cite{Pletser}. We prove the following theorem using only elementary methods, which implies there are only trivial solutions to the equation where $y$ is a prime power. In other words, there are no non-trivial solutions unless $y$ is divisible by at least two different primes.
\begin{theorem}
For $x,k,r \in\mathbb{Z^+}$ and $p$ prime, all solutions to the equation
\begin{equation}
x^3 + (x+1)^3 + ... + (x+k-1)^3 = (p^r)^2
\label{primeeq}
\end{equation}

are $(x,k,p,r) = (p^{2c},1,p,3c), (1,2,3,1)$ where $c$ is an arbitrary positive integer. \label{theorem1}
\end{theorem}
Motivation for proof is the Pythagorean triple appearing after a few steps on Eq. (\ref{maineq}),
$$
\sum_{i = x}^{x+k-1} i^3 = \sum_{i = 1}^{x+k-1} i^3 - \sum_{i = 1}^{x-1} i^3 = \Delta_{x+k-1}^2 - \Delta_{x-1}^2= y^2
$$
\begin{equation}
\bigg[\frac{(x+k-1)(x+k)}{2}\bigg]^2 = p^{2r} + \bigg[\frac{(x-1)x}{2}\bigg]^2
\label{mainpythagoreanprimeeq}
\end{equation}
In cooperation with Pythagorean triples, another key element is the $v_{p}(x)$ notation, denoting the non-negative integer $r$ such that $p^r\mid x$ but $p^{r+1}\nmid x$. However, including $v_{p}(\Delta_{x-1}^2)$ in solution will not be valid (Sections \ref{evencase} and \ref{oddcase}) since the expression becomes undefined when $x=1$. Therefore, a separate investigation of $x=1$ is imperative.
\subsection{For $x$ equals one}
Rewriting (\ref{mainpythagoreanprimeeq}) with the assumption $x=1$ gives,
$$
\bigg[\frac{k(k+1)}{2}\bigg]^2 = p^{2r} \Longleftrightarrow \frac{k(k+1)}{2} = p^r.
$$
Last equation implies either $k = p^a = 2p^{r-a}-1$ or $k = 2p^a = p^{r-a}-1$. In both cases, modulo $p$ yields $a = 0$ and thereby, we obtain $(x,k,y) = (1,1,1), (1,2,3)$. In this case, the only valid solution for $y$ to be a prime power is $(x,k,p,r) = (1,2,3,1)$.

\subsection{For $x$ bigger than one}\label{xbiggerthanone}
In Eq. (\ref{primeeq}), the terms can be grouped as shown,
$$
\big[x^3 + (x+k-1)^3\big] + \big[(x+1)^3 + (x+k-2)^3\big] + \big[(x+2)^3 + (x+k-3)^3\big] + ... = (p^r)^2 .
$$
For $k$ even, each square bracket will be divisible by $2x+k-1$, hence $2x+k-1\mid (p^r)^2$ and $2x+k-1 = p^t$ for $0<t<2r$.\\
For $k$ odd, square brackets will be divisible by $2x+k-1$ again, but also a median term exists, which is $\big(\frac{2x+k-1}{2}\big)^3$. Hence $\frac{2x+k-1}{2}$ will divide right side. Therefore, $2x+k-1 = 2p^t$ for $0<t<2r$.\\
These cases should be considered separately by parity of $k$.
\begin{remark}
Previous arguments are true, independent of the power two on the right side, as long as the common power of terms (Which is three in this case) preserves its parity.
\end{remark}
\section{For $k$ even} \label{evencase}
We assume $k = 2l$ for $l\in \mathbb{Z^+}$. The result above gives $2x+2l-1 \mid (p^r)^2$ and underlines the fact that
\begin{equation}
    2x+2l-1 = p^t\ , \ 0<t<2r.
    \label{primepower1}
\end{equation}
Rewrite Eq. (\ref{mainpythagoreanprimeeq}) and apply (\ref{primepower1}), it yields,
\begin{equation}
 \left[\frac{(p^t-x)(p^t-(x-1))}{2} \right]^2 = (p^r)^2 + \left[\frac{(x-1)x}{2} \right]^2.
 \label{kevenmaineq}
\end{equation}
We aim to take advantage of several facts, which are formulated in Eq. (\ref{qlessteven}) and (\ref{triangleequiivkeven}):
\begin{equation}
\ q = v_{p}\Bigg(\frac{(x-1)x}{2}\Bigg), q < t. \label{qlessteven}
\end{equation}
\begin{equation}
    \frac{(p^t-x)(p^t-(x-1))}{2} \equiv \frac{(x-1)x}{2} \Mod{p^t}.
    \label{triangleequiivkeven}
\end{equation}
Eq. (\ref{qlessteven}) is derived from the following inequalities whereas (\ref{triangleequiivkeven}) is the explicit result of a simple $\Mod{p^t}$ investigation of Eq. (\ref{kevenmaineq}). 
$$
v_{p}\Bigg(\frac{(x-1)x}{2}\Bigg) \leq \max\{v_{p}(x-1), v_{p}(x)\} < v_{p}(2x+2l-1) = v_{p}(p^t) = t.
$$

Combining (\ref{qlessteven}) and (\ref{triangleequiivkeven}) we have
\begin{equation}
    q = v_{p}\bigg(\frac{(p^t-x)(p^t-(x-1))}{2}\bigg) = v_{p}\bigg(\frac{(x-1)x}{2} \bigg).
    \label{vpequalkeven}
\end{equation}
In Eq. (\ref{kevenmaineq}) with modulo $p^{2q}$, we see $p^{2q}\mid p^{2r}$. Let Eq. (\ref{kevenmaineq}) become $a^2 = (p^{r-q})^2 + b^2$ when every term is divided by $p^{2q}$. From (\ref{vpequalkeven}) $a$ and $b$ has no common divisor greater than $1$, hence $(a,p^{r-q},b)$ is a primitive Pythagorean triple. Thus, we obtain two cases by identifying general forms them as follows.
\begin{alignat*}{3}
    &a = m^2+n^2 &\ \ \ \ \ \ &a = m^2+n^2\\
    &b = m^2-n^2 &\ \ \ \ \ \ &b = 2mn\\
    &p^{r-q} = 2mn &\ \ \ \ \ \ &p^{r-q} = m^2-n^2
\end{alignat*}
\subsection{First Case}\label{pythagoreanfirstcasefirst}
The equation $p^{r-q} = 2mn$ yields $p = 2$, and therefore $m$ and $n$ are perfect powers of two. However, since $a$ and $b$ are coprime, only solutions are $n=1$ and $a-b = 2$. From there 
$$
\frac{(p^t-x)(p^t-(x-1))}{2} - \frac{(x-1)x}{2} = (a-b)p^q = 2p^q = 2^{q+1}
$$
and we also have
$$
 \frac{(p^t-x)(p^t-(x-1))}{2} - \frac{(x-1)x}{2} = \frac{p^{t}(p^{t} - 2x + 1)}{2} = 2^{q+1}
 $$
Applying (\ref{primepower1}) it follows that,
$$
\underline{2^{q+1}} = \frac{p^{t}(p^{t} - 2x + 1)}{2} =\frac{p^{t}(2x+2l-1-2x+1)}{2} =p^{t}l = \underline{2^{t}l}
 $$
From underlined terms in the last equation and Eq. (\ref{qlessteven}), we have $t = q+1$ and $l=1$. Rewriting Eq. (\ref{primepower1}) with $l=1$ gives $2x+1 = 2^t$. This implies $t=0$ by parity and contradicts the range of $t$ at Eq. (\ref{primepower1}) and $x>1$.

\subsection{Second Case}\label{pythagoreansecondcasefirst}
For the second system, we presumed that $p^{r-q}$ is the product of $m-n$ and $m+n$. There exists $c$ and $d$ such that $0\leq c<d$, $m-n = p^c$ and $m+n = p^d$. Then, the following equations hold.
$$
 \frac{p^{2c} + p^{2d}}{2} = \frac{(m-n)^2 + (m+n)^2}{2} = m^2+n^2 = a
 $$
Notice that any positive value of $c$ causes the leftmost side of the equation to give $0$ in $\Mod{p}$. Hence, the rightmost side of the equation will give $0$ in $\Mod{p}$ which contradicts with $(a,p)=1$. Consequently, $c$ must be zero. The equation $a-b = (m-n)^2 = p^c$ then becomes $a-b=1$. It follows that 
$$
 \frac{(p^t-x)(p^t-(x-1))}{2} - \frac{(x-1)x}{2} = \frac{p^{t}(p^{t} - 2x + 1)}{2} = (a-b)p^q = p^q.
$$
Similar to the previous case, we apply (\ref{primepower1}) 
$$
\frac{p^{t}(p^{t} - 2x + 1)}{2} =\frac{p^{t}(2x+2l-1-2x+1)}{2} = p^{t}l =  p^q.
$$
Since $l$ is a positive integer, $q\geq t$ must hold, which contradicts with Eq. (\ref{qlessteven}). In conclusion, there are no other solutions to Eq. (\ref{primeeq}) when $k$ is even and $x>1$.

\section{For $k$ odd \label{oddcase}}
In this case, we assume $k=2l-1$ for $l \in \mathbb{Z^+}$. In Section (\ref{xbiggerthanone}), it is stated that $2x+k-1\mid 2p^{2r}$ and therefore
\begin{equation}
x+l-1 = p^t,\ 0<t<2r
\label{primepower2}
\end{equation}
The solution follows similar steps to Section \ref{evencase}. Using Eq. (\ref{primepower2}) in Eq. (\ref{mainpythagoreanprimeeq}) gives
\begin{equation}
 \left[\frac{(2p^t - x)(2p^t - (x-1))}{2} \right]^2 = (p^r)^2 + \left[\frac{(x-1)x}{2} \right]^2
 \label{koddmaineq}
\end{equation}
$$
\underline{q} = v_{p}\bigg(\frac{(x-1)x}{2} \bigg) \leq \max\{v_{p}(x-1),v_{p}(x)\} \leq v_{p}(x+l-1) = v_{p}(p^t) = \underline{t}.
$$
The last row shows $q\leq t$. The cases $q<t$ and $q=t$ shall be investigated separately.
\subsection{For $q$ less than $t$}\label{chapterqlessthant}
Combining $q<t$ and $\Delta_{x+k-1} \equiv \Delta_{x-1} \Mod{p^t}$ from Eq. (\ref{koddmaineq}) yields the following. 
\begin{equation}
    q = v_{p}\bigg( \frac{(2p^t - x)(2p^t - (x-1))}{2}\bigg) = v_{p}\bigg(\frac{(x-1)x}{2}\bigg)
    \label{vpequalkodd}
\end{equation}
In Eq. (\ref{koddmaineq}), $p^{2q}$ divides $p^{2r}$ and let irreducible equation be $a^2 = (p^{r-q})^2 + b^2$. It is a Pythagorean triple and also primitive because of (\ref{vpequalkodd}). By identifying general forms of these triples, we obtain two cases that are analogous to the ones in Section (2).
\subsubsection{First Case}\label{pythagoreanfirstcasesecond}
The same procedure applied in Section (\ref{pythagoreanfirstcasefirst}) easily yields $p=2$, $a-b=2$ and 
$$
(a-b)p^q = \underline{2^{q+1}} = \frac{(2p^t - x)(2p^t - (x-1))}{2} - \frac{(x-1)x}{2} =  p^t(2p^t-2x+1) = \underline{p^t(2l-1)}.
$$
Since we assumed $q<t$, the underlined parts of the equation hold if and only if $k=l=1$. Plugging these in Eq. (\ref{primepower2}) yields $t=2c$, $r=3c$ for an arbitrary positive integer $c$ and the solution $(x,k,p,r)$ becomes $(2^{2c},1,2,3c)$.
\subsubsection{Second Case}
Similar to Section (\ref{pythagoreansecondcasefirst}), we obtain $a-b = 1$ and thus,
$$
\frac{(2p^t - x)(2p^t - (x-1))}{2} - \frac{(x-1)x}{2} = (a-b)p^q = p^q.
$$
Equation simplifies to $p^t(2p^t-2x+1) = p^q$ which contradicts the assumption $q<t$ in Section (\ref{chapterqlessthant}). Therefore there is no solution in this case.
\subsection{For $q$ equals $t$}
If $q=t$, then $p^t$ divides every term in Eq. (\ref{koddmaineq}). From (\ref{primepower2}) we have $p^t<2p^t - (x-1) < 2p^t$ and thereby $p^t\nmid (2p^t - (x-1))$. This implies $p^t\mid (2p^t -x)$ and according to (\ref{primepower2}), condition holds only when $x=p^t$ and $l=1$. Eq. (\ref{koddmaineq}) becomes
$$
 \left[ \frac{p^t(p^t+1)}{2}\right]^2 = (p^r)^2 + \left[ \frac{p^t(p^t-1)}{2}\right]^2
$$
If $p=2$, then clearly $q<t$ contradicts. Hence, $p>2$, and every term is divisible by $p^t$. The above equation simplifies to
$$
(p^{r-t})^2 = \left[\frac{p^t+1}{2} \right]^2 - \left[\frac{p^t-1}{2} \right]^2 = p^t.
$$
We get $2(r-t) = t$, $2r = 3t$. For an arbitrary positive integer $c$ and $p>2$ it yields the solution $(x,k,p,r) = (p^{2c},1,p,3c)$. The case $p=2$ is obtained from the Section (\ref{pythagoreanfirstcasesecond}). Thus, all solutions are as stated in Theorem \ref{theorem1} $\blacksquare$.

\paragraph*{Acknowledgements.} I especially thank Prof. Dr. Gökhan Soydan for their invaluable suggestions and mentoring about my study. I would like to thank Professor Nikos Tzanakis and Dr. Paul Voutier for sharing their translation of `Ein Satz Uber Die Diophantische Gleichung  $ Ax^2 - By^4 = C$ , ($C = 1, 2, 4$)' to English. Their translation was very useful for me.

\end{document}